\documentclass[12pt]{article}
\usepackage{graphics,amsmath,amssymb}
\usepackage{latexsym}
\usepackage{epsfig}
\usepackage{hyperref}
\usepackage{fancybox}
\usepackage{xcolor}

\def\epsilon{\varepsilon}

\def\phi{\varphi}

\newtheorem{theorem}{Theorem}[section]
\newtheorem{lemma}[theorem]{Lemma}

\newtheorem{definition}[theorem]{Definition}

\newtheorem{remark}[theorem]{Remark}

\def\N{{\mathbb N}}
\def\C{{\mathbb C}}
\def\R{{\mathbb R}}

\def\longto#1{%
\mathrel{\mathop{\kern0pt\longrightarrow}\limits_{#1}}}

\newenvironment{Proof}{\removelastskip\par\medskip
\noindent{\em Proof.} \rm}{\penalty-20\null$\square$\par\medbreak}

\title{\bf   Weak solutions for time-fractional  
\\
evolution equations in Hilbert spaces}

\author{Paola Loreti
\thanks{Dipartimento di Scienze di Base e Applicate per l'Ingegneria,
Sapienza Universit\`a di Roma,
Via Antonio Scarpa 16, 00161 Roma (Italy); e-mail: 
$<$paola.loreti@uniroma1.it$>$ }
\and Daniela Sforza
\thanks{Dipartimento di Scienze di Base e Applicate per l'Ingegneria, 
Sapienza Universit\`a di Roma,
Via Antonio Scarpa 16, 00161 Roma (Italy); e-mail: 
$<$daniela.sforza@uniroma1.it$>$ }}

\begin{document}
\date{}

\maketitle

\begin{abstract}
We introduce a notion of  weak solution for abstract fractional differential equations, motivated by the definition of Caputo derivative. 
We prove existence results for  weak and strong solutions. We also give two examples as application of our results: time-fractional wave equations and time-fractional Petrovsky systems. 

\end{abstract}



\bigskip
\noindent

\section{Introduction}

%
%
%
%
%

Our aim  is to establish existence and regularity results regarding weak solutions  of
the abstract fractional  equation 
\begin{equation}\label{eq:weakpI}
\partial_t^{\alpha}u+A u=0
\qquad \mbox{in}\ (0,T),
\end{equation}
where $u$ takes its values in a Hilbert space $H$, $A$ is a densely defined, linear self-adjoint positive operator on $H$ and the Caputo derivative $\partial_t^{\alpha}u$ of order $\alpha\in(1,2)$ is defined as
 \begin{equation*}
\partial_t^{\alpha}u(t)
=\frac1{\Gamma(2-\alpha)}\int_0^t (t-\tau)^{1-\alpha}\frac{d^2u}{d\tau^2}(\tau)\ d\tau\,.
\end{equation*} 
Our study is inspired from the literature concernig partial differential equations of the type \ref{eq:weakpI} in the case the operator $A$ is a concrete differential operator. Typically, the operator $A$ is the Laplacian or, more general, a uniformly elliptic operator. To our knowledge, the biharmonic operator has not been studied yet in the framework of the time-fractional differential equations.   

We mention \cite{Fujita}  where the author gives a first analysis about evolution equations, related to fractional differential ones, that interpolates the heat equation and the wave equation in the domain $(0,+\infty)\times\R$.    
In the survey \cite{MG} the interested reader can find a historical note on the origin of the fractional derivative in the Caputo sense. Morever, the authors provide some arguments showing the usefulness of the Caputo derivatives in the theory of viscoelasticity. 

To begin with, we introduce strong solutions  of \eqref{eq:weakpI} as functions $u$ belonging to $C([0,T];D(A))\cap C^1([0,T];H)$, $\partial_t^{\alpha}u\in C([0,T];H)$ and $u$ satisfies
equation \eqref{eq:weakpI} for any $t\in [0,T]$.

In this paper we are interested to solutions with lower regularity, so the main difficulty
is an appropriate definition for  weak solutions. The definition of weak solutions that we propose is suggested by the formula 
\begin{equation}\label{eq:}
\partial_t^{\alpha}u(t)=\frac{d}{dt}I^{2-\alpha}\big(u'-u'(0)\big)(t)
\,,
\end{equation}
where $I^{2-\alpha}$ is the Riemann--Liouville operator of order $2-\alpha$, see formula \eqref{eq:RLfi} in Subsection \ref{s:frac-d}.

We define a weak solution of the fractional differential equation \eqref{eq:weakpI} as a
function $u$ belonging to $C([0,T];D(\sqrt{A}))$, $u'\in L^2(0,T;H)\cap C([0,T];D(A^{-\theta}))$,  for some $\theta\in(0,1)$, 
and for any $v\in D(\sqrt{A})$ one has $\langle I^{2-\alpha}\big(u'-u'(0)\big)(t),v\rangle\in C^1([0,T])$ and 
\begin{equation}\label{eq:w-int-I}
\frac{d}{dt}\langle I^{2-\alpha}\big(u'-u'(0)\big)(t),v\rangle+\langle \sqrt{A}u(t),\sqrt{A} v \rangle=0
\qquad t\in [0,T]
\,.
\end{equation}
Our main result is the following (see Section \ref{s:reg})
\begin{theorem}\label{th:reg-l2}
If $u_0\in D(\sqrt{A})$ and $u_1\in H$, then the function
\begin{equation}\label{eq:def-u0-I}
u(t)=\sum_{n=1}^\infty\big[  \langle u_0,e_n\rangle E_{\alpha,1}(-\lambda_nt^\alpha) 
+\langle u_1,e_n\rangle t E_{\alpha,2}(-\lambda_nt^\alpha)\big]e_n
\end{equation}
is the unique weak solution of \eqref{eq:weakpI} satisfying
 the initial conditions
\begin{equation*}
u(0)=u_{0},\quad
u'(0)=u_{1}.
\end{equation*}
\end{theorem}
Formula \eqref{eq:def-u0-I} is suggested by the spectral approach given in \cite{SY,KY}, where the authors deal with uniformly elliptic operators in a different setting. It is worthwile to mention the paper \cite{HMP} (see also references therein), where  the authors give the representation of classical solutions by means of the $\alpha$-resolvent family.

Our abstract results can be applied to various partial differential equations to obtain the existence of weak solutions.
The last section is devoted to the discussion of two examples: time-fractional wave equations and time-fractional Petrovsky systems.

\section{Preliminaries}\label{s:pre}
In this section we collect some notations, definitions and known results that we use to prove our main results.

\subsection{Abstract operators}\label{s:}
Let $H$ be a real Hilbert space with inner product
$\langle \cdot \, ,\, \cdot \rangle$ and norm $\| \cdot \|$.
$A$ is a linear self-adjoint positive
  operator on $H$ with dense domain $D(A)$, satisfying
\begin{equation}\label{eq:opA}
\langle Ax,x\rangle\ge a\|x\|^2\qquad\forall x\in D(A)
\end{equation}
for some $a>0$. We assume that the spectrum of $A$ consists  of a sequence of positive eigenvalues $\{\lambda_n\}_{n\in\N}$ such that $\lambda_n\longto{n}\infty$. Moreover, $\lambda_n$ are all distinct numbers, whence the eigenspace generated by every single $\lambda_n$  has dimension one. 

Moreover, the eigenfunctions $e_n$ of $A$  \big($A e_n=\lambda_n e_n$\big) constitute an orthonormal basis of $H$.

The fractional powers $A^{\theta}$ are defined for $\theta>0$, see e.g. \cite{Pazy,Lunardi}. 
The domain $D(A^{\theta})$ of $A^{\theta}$ consists of  $u\in H$ such that
\begin{equation*}
\sum_{n=1}^\infty \lambda_n^{2\theta} |\langle u,e_n\rangle|^2<+\infty
\end{equation*}
and
\begin{equation*}
A^{\theta} u=\sum_{n=1}^\infty \lambda_n^{\theta} \langle u,e_n\rangle e_n,
\quad
u\in D(A^{\theta}).
\end{equation*}
Moreover $D(A^{\theta})$ is a Hilbert space with the norm given by
\begin{equation}\label{eq:norm-frac}
\|u\|_{D(A^{\theta})}=\|A^{\theta} u\|=\left(\sum_{n=1}^\infty \lambda_n^{2\theta} |\langle u,e_n\rangle|^2\right)^{1/2}
\quad
u\in D(A^{\theta}),
\end{equation}
and for any $0<\theta_1<\theta_2$ we have $D(A^{\theta_2})\subset D(A^{\theta_1})$.

In particular,  the  norm of the space $D(\sqrt{A})$   is given by
\begin{equation}\label{eq:H^2_0}
\|u\|_{D(\sqrt{A})}=\|\sqrt{A} u\|=\left(\sum_{n=1}^\infty \lambda_n |\langle u,e_n\rangle|^2\right)^{1/2}
\qquad
u\in D(\sqrt{A}).
\end{equation}
If we identify the dual $H'$ with $H$ itself, then we have $D(A^{\theta})\subset H\subset(D(A^{\theta}))'$.
From now on we set 
\begin{equation}
D(A^{-\theta}):=(D(A^{\theta}))',
\end{equation}
whose elements are bounded linear functionals on $D(A^{\theta})$. If $u\in D(A^{-\theta})$ and $\varphi\in D(A^{\theta})$ the value $u(\varphi)$ is denoted  
by 
\begin{equation}\label{eq:duality}
\langle u,\varphi\rangle_{-\theta,\theta}:=u(\varphi)\,.
\end{equation}
 In addition, 
 $D(A^{-\theta})$ is a Hilbert space with the norm given by
\begin{equation}\label{eq:norm-theta}
\|u\|_{D(A^{-\theta})}=\left(\sum_{n=1}^\infty \lambda_n^{-2\theta} |\langle u,e_n\rangle_{-\theta,\theta}|^2\right)^{1/2}
\quad
u\in D(A^{-\theta})
\,,
\end{equation}
and for any $0<\theta_1<\theta_2$ we have $D(A^{-\theta_1})\subset D(A^{-\theta_2})$.
We also recall that 
\begin{equation}\label{eq:-theta}
\langle u,\varphi\rangle_{-\theta,\theta}=\langle u,\varphi\rangle
\qquad \mbox{for}\ u\in H\,,\varphi\in D(A^{\theta}),
\end{equation}
see e.g. \cite[Chapitre V]{B}.
%

\subsection{Fractional derivatives}\label{s:frac-d}
\begin{definition}
For any $\beta>0$
we denote the Riemann--Liouville fractional integral operator of order $\beta$ by
\begin{equation}\label{eq:RLfi}
I^{\beta}(f)(t)=\frac1{\Gamma(\beta)}\int_0^t (t-\tau)^{\beta-1}f(\tau)\ d\tau, 
\qquad
f\in L^1(0,T),
\ \mbox{a.e.}\ t\in(0,T), 
\end{equation}
where $T>0$ and $\Gamma (\beta)=\int_0^\infty t^{\beta-1}e^{-t}\ dt$ is the Euler Gamma function. 
\end{definition}
The Caputo fractional derivative of order $\alpha\in(1,2)$ is given by
\begin{equation}\label{eq:der-frac}
\partial_t^{\alpha}f(t)
=\frac1{\Gamma(2-\alpha)}\int_0^t (t-\tau)^{1-\alpha}\frac{d^2f}{d\tau^2}(\tau)\ d\tau\,.
\end{equation}
By means of the Riemann--Liouville  integral operator $I^{2-\alpha}$ we can write
\begin{equation}\label{eq:der-fracI}
\partial_t^{\alpha}f(t)
=
I^{2-\alpha}\Big(\frac{d^2f}{dt^2}\Big)(t)
\,.
\end{equation}
We also note that if $f'$ is absolutely continuous, then
\begin{equation}\label{eq:der-fracI1}
\partial_t^{\alpha}f(t)=\frac{d}{dt}I^{2-\alpha}\big(f'-f'(0)\big)(t)
\,.
\end{equation}
For arbitrary constants $\alpha,\beta> 0$, we denote the Mittag--Leffler functions by
\begin{equation}
E_{\alpha,\beta}(z):= \sum_{k=0}^\infty\frac{z^k}{\Gamma(\alpha k+\beta)},
\quad z\in\C.
\end{equation}
By the power series, one can note that $E_{\alpha,\beta}(z)$ is an entire function
of $z\in\C$. We note that $E_{\alpha,1}(0)=1$.

The proof of the following result can be found in \cite[p. 35]{Pod}, see also \cite[Lemma 3.1]{SY}.
In the following we denote the Laplace transform of a function $f(t)$ by the symbol
\begin{equation}\label{eq:laplace}
{\cal L}[f(t)](z):=\int_0^{\infty}e^{-zt}f(t)\ dt
\qquad z\in\C.
\end{equation}
\begin{lemma}
\begin{enumerate}
\item 
Let $\alpha\in(1,2)$ and $\beta>0$ be. Then for any $\mu\in\R$ such that $\pi\alpha/2<\mu<\pi$ there exists a constant $C=C(\alpha,\beta,\mu)>0$ such that 
\begin{equation}\label{eq:stimeE}
\big|E_{\alpha,\beta}(z)\big|\le \frac{C}{1+|z|},
\qquad z\in\C,\ \mu\le|\arg(z)|\le\pi.
\end{equation}
\item
For $\alpha\,,\beta\,,\lambda>0$ one has
\begin{equation}\label{eq:LTML}
{\cal L}\big[t^{\beta-1}E_{\alpha,\beta}(-\lambda t^\alpha)\big)\big](z)=\frac{z^{\alpha-\beta}}{z^\alpha+\lambda}
\qquad \Re z>\lambda^{\frac1\alpha}\,.
\end{equation}
\item 
If $\alpha\,,\lambda>0$, then we have
\begin{equation}\label{eq:Ea1}
\frac{d}{dt}E_{\alpha,1}(-\lambda t^{\alpha})=-\lambda t^{\alpha-1}E_{\alpha,\alpha}(-\lambda t^{\alpha}), 
\qquad t>0,
\end{equation}
\begin{equation}\label{eq:Eaa1}
\frac{d}{dt}\Big(t^{k}E_{\alpha,k+1}(-\lambda t^{\alpha})\Big)=t^{k-1}E_{\alpha,k}(-\lambda t^{\alpha}), 
\qquad k\in\N\,,t\ge0,
\end{equation}
\begin{equation}\label{eq:Eaaa}
\frac{d}{dt}\Big(t^{\alpha-1}E_{\alpha,\alpha}(-\lambda t^{\alpha})\Big)=t^{\alpha-2}E_{\alpha,\alpha-1}(-\lambda t^{\alpha}), 
\qquad t\ge0.
\end{equation}
\end{enumerate}
\end{lemma}

We also exhibit an elementary result that will be useful in the estimates.
\begin{lemma} For any $0<\beta<1$ the function $x\to\frac{x^\beta}{1+x}$ gains its maximum on $[0,+\infty[$ at point $\frac\beta{1-\beta}$ and the maximum value is given by
\begin{equation}\label{eq:maxbeta}
\max_{x\ge0}\frac{x^\beta}{1+x}=\beta^\beta(1-\beta)^{1-\beta},
\qquad\beta\in (0,1)
\,.
\end{equation}

\end{lemma}
Now we recall the definition of fractional vector-valued Sobolev spaces.
For $\beta\in (0,1)$, $T>0$ and a Hilbert space $H$, endowed with the norm $\|\cdot\|_H$, $H^\beta(0,T;H)$ is the space of all $u\in L^2(0,T;H)$ such that 
\begin{equation}\label{eq:gagliardo}
[u]_{H^\beta(0,T;H)}:=\left(\int_0^T\int_0^T\frac{\|u(t)-u(\tau)\|_H^2}{|t-\tau|^{1+2\beta}}\ dtd\tau\right)^{1/2}<+\infty\,,
\end{equation}
that is $[u]_{H^\beta(0,T;H)}$ is the so-called Gagliardo semi-norm of $u$.
$H^\beta(0,T;H)$ is endowed with the norm
\begin{equation}\label{eq:defHs}
\|\cdot\|_{H^\beta(0,T;H)}:=\|\cdot\|_{L^2(0,T;H)}+[\ \cdot\ ]_{H^\beta(0,T;H)}.
\end{equation}
The following result is a generalization  to the case of vector valued functions of  
\cite[Theorem 2.1]{GLY}. 
We will use the symbol $\sim$ to indicate equivalent norms.
\begin{theorem}\label{th:R-Lop}
Let $H$ be a separable Hilbert space.
\begin{itemize}
\item[(i)] 
The Riemann--Liouville operator $I^{\beta}:L^2(0,T;H)\to L^2(0,T;H)$, $0<\beta\le1$, is injective and the range ${\cal R}(I^{\beta})$ of $I^{\beta}$ is given by
\begin{equation}\label{eq:range}
{\cal R}(I^{\beta})=
\begin{cases}H^\beta(0,T;H), \hskip5.5cm 0<\beta<\frac12,
\\
\Big\{v\in H^{\frac12}(0,T;H): \int_0^Tt^{-1}|v(t)|^2 dt<\infty\Big\},
\qquad \beta=\frac12,
\\
_0H^\beta(0,T;H), \hskip5.2cm \frac12<\beta\le1,
\end{cases}
\end{equation}
where $_0H^\beta(0,T)=\{u\in H^\beta(0,T): u(0)=0\}$.
\item[(ii)] 
For the Riemann--Liouville operator $I^{\beta}$ and its inverse operator $I^{-\beta}$ the norm equivalences 
\begin{equation}\label{eq:R-Lop}
\begin{split}
\|I^{\beta}(u)\|_{H^\beta(0,T;H)}
&\sim\|u\|_{L^2(0,T;H)},
\qquad u\in L^2(0,T;H),
\\
\|I^{-\beta}(v)\|_{L^2(0,T;H)}
&\sim\|v\|_{H^\beta(0,T;H)},
\qquad v\in {\cal R}(I^{\beta}),
\end{split}
\end{equation}
hold true.
\end{itemize}
\end{theorem}

\section{Existence and regularity of solutions}\label{s:reg}

First, we introduce  the notions of weak and strong solutions,  see \cite[]{}. Let $\alpha\in(1,2)$ and $T>0$.
\begin{definition}\label{de:wss}
\begin{enumerate}
\item[\bf 1.] 
A function $u$ is called a weak solution of the abstract time-fractional equation
\begin{equation}\label{eq:weakp}
\partial_t^{\alpha}u+A u=0
\end{equation}
if  $u\in C([0,T];D(\sqrt{A}))$, $u'\in L^2(0,T;H)\cap C([0,T];D(A^{-\theta}))$,  for some $\theta\in(0,1)$, 
and for any $v\in D(\sqrt{A})$ one has $\langle I^{2-\alpha}\big(u'-u'(0)\big)(t),v\rangle\in C^1([0,T])$ and 
\begin{equation}\label{eq:w-int}
\frac{d}{dt}\langle I^{2-\alpha}\big(u'-u'(0)\big)(t),v\rangle+\langle \sqrt{A}u(t),\sqrt{A} v \rangle=0
\qquad t\in [0,T]
\,.
\end{equation}
\item [\bf 2.]
A function $u$ is called a strong solution of \eqref{eq:weakp} if $u\in C([0,T];D(A))\cap C^1([0,T];H)$, $\partial_t^{\alpha}u\in C([0,T];H)$ and satisfies
equation \eqref{eq:weakp} for any $t\in [0,T]$.
\end{enumerate}
\end{definition}
\begin{remark}
We observe that a strong solution is also a weak solution.

Also, we note that for  a weak solution $u$ of \eqref{eq:weakp} we have $\partial_t^{\beta}u\in H^{1-\beta}(0,T;H)$, $\beta\in(0,1)$, where
the Caputo fractional derivative of order $\beta\in(0,1)$ is defined as
\begin{equation}\label{eq:der-frac01}
\partial_t^{\beta}u(t)
=\frac1{\Gamma(1-\beta)}\int_0^t (t-\tau)^{-\beta}u'(\tau)\ d\tau
=I^{1-\beta}\big(u'\big)(t)\,.
\end{equation}
Indeed, since $u'\in L^2(0,T;H)$ we can apply Theorem \ref{th:R-Lop} to obtain $I^{1-\beta}\big(u'\big)\in H^{1-\beta}(0,T;H)$, that is
$\partial_t^{\beta}u\in H^{1-\beta}(0,T;H)$.

In particular, for $\beta=\alpha/2$ and $\beta=1-\alpha/2$ we have $\partial_t^{\alpha/2}u\in H^{1-\alpha/2}(0,T;H)$ and 
$\partial_t^{1-\alpha/2}u\in H^{\alpha/2}(0,T;H)$ respectively.
\end{remark}
A classical approach to solve scalar fractional differential equations is using Laplace transform methods. 
\begin{lemma}\label{le:x-y}
For any $\lambda>0$ and $x,y\in\R$ the solution of  problem
\begin{equation}\label{Aeq:cauchy100}
\begin{cases}
\displaystyle
\partial_t^{\alpha}u(t)+\lambda u(t)=0\,,
\quad t\ge0,
\\
u(0)=x,\quad
u'(0)=y,
\end{cases}
\end{equation}
is given by
\begin{equation}\label{eq:x-y}
u(t)=x\ E_{\alpha,1}(-\lambda t^\alpha) 
+y\ t E_{\alpha,2}(-\lambda t^\alpha) ,
\qquad  t\ge0
\,.
\end{equation}
\end{lemma}
\begin{Proof}
First, we compute the Laplace transform of $\partial_t^{\alpha}u(t)$, that is 
\begin{multline*}
{\cal L}[\partial_t^{\alpha}u(t)](z)=z^{\alpha-2}{\cal L}\big[u''(t)\big](z)
=z^{\alpha-2}\big(z^2{\cal L}[u(t)](z)-x\ z-y\big)
\\
=z^{\alpha}{\cal L}[u(t)](z)-x\ z^{\alpha-1}-y\ z^{\alpha-2}
\,.\end{multline*}
Therefore, taking the Laplace transform of the equation in \eqref{Aeq:cauchy100} we obtain
\begin{equation*}
{\cal L}[\partial_t^{\alpha}u(t)](z)+\lambda{\cal L}[u(t)](z)
=(z^{\alpha}+\lambda){\cal L}[u(t)](z)-x\ z^{\alpha-1}-y\ z^{\alpha-2}
=0
\end{equation*}
and hence
\begin{equation*}
{\cal L}[u(t)](z)=x\ \frac{z^{\alpha-1}}{z^{\alpha}+\lambda}+y\ \frac{z^{\alpha-2}}{z^{\alpha}+\lambda}
\,.
\end{equation*}
From \eqref{eq:LTML} we have
\begin{equation*}
{\cal L}\big[E_{\alpha,1}(-\lambda t^\alpha)\big](z)=\frac{z^{\alpha-1}}{z^{\alpha}+\lambda},
\qquad
{\cal L}\big[t E_{\alpha,2}(-\lambda t^\alpha)\big](z)=\frac{z^{\alpha-2}}{z^{\alpha}+\lambda},
\quad \Re z>\lambda^{1/\alpha},
\end{equation*}
whence we deduce
\begin{multline*}
{\cal L}[u(t)](z)
=x\ {\cal L}\big[E_{\alpha,1}(-\lambda t^\alpha)\big](z)+y\ {\cal L}\big[t E_{\alpha,2}(-\lambda t^\alpha)\big](z)
\\
={\cal L}\big[x\ E_{\alpha,1}(-\lambda t^\alpha)+y\ t E_{\alpha,2}(-\lambda t^\alpha)\big](z)
\,.
\end{multline*}
In conclusion, by Lerch's Theorem that assures the uniqueness of inverse Laplace transforms $u$ is given by \eqref{eq:x-y}.
\end{Proof}

\begin{theorem}\label{th:reg-l2}
\begin{enumerate}
\item[(i)] 
If $u_0\in D(\sqrt{A})$ and $u_1\in H$, then the function
\begin{equation}\label{eq:def-u0}
u(t)=\sum_{n=1}^\infty\big[  \langle u_0,e_n\rangle E_{\alpha,1}(-\lambda_nt^\alpha) 
+\langle u_1,e_n\rangle t E_{\alpha,2}(-\lambda_nt^\alpha)\big]e_n
\end{equation}
is the unique weak solution of \eqref{eq:weakp} satisfying
 the initial conditions
\begin{equation}\label{eq:ini-cond}
u(0)=u_{0},\quad
u_t(0)=u_{1}.
\end{equation}
In addition
\begin{equation}\label{eq:def-u_t}
u'(t)=\sum_{n=1}^\infty\big[-\lambda_n \langle u_0,e_n\rangle t^{\alpha-1}E_{\alpha,\alpha}(-\lambda_nt^\alpha)
+\langle u_1,e_n\rangle E_{\alpha,1}(-\lambda_nt^\alpha)\big]e_n
\,,
\end{equation}
and $u'\in C([0,T];D(A^{-\theta}))$ for  $\theta\in\big(\frac{2-\alpha}{2\alpha},\frac12\big)$.
\item[(ii)]
For $u_0\in D({A})$ and $u_1\in D(\sqrt{A})$
the weak solution given by \eqref{eq:def-u0} is a strong one and 
\begin{equation}\label{eq:def-u-alpha}
\partial_t^{\alpha}u(t)=-\sum_{n=1}^\infty\big[\lambda_n \langle u_0,e_n\rangle E_{\alpha,1}(-\lambda_nt^\alpha)
+\lambda_n\langle u_1,e_n\rangle t E_{\alpha,2}(-\lambda_nt^\alpha)\big]e_n
\,.
\end{equation}
\end{enumerate}

\end{theorem}
\begin{Proof}
To begin with, we intend to justify the expression of the solution given by the series \eqref{eq:def-u0}.
To this end we search the solution in the form 
\begin{equation*}
u(t)=\sum_{n=1}^\infty\ u_n(t)e_n
\end{equation*}
where the functions $u_n(t)=\langle u(t),e_n\rangle$ are unknown. First, we assume that there exists a strong solution $u$, see Definition \ref{de:wss}-(2).  By means of the scalar product in $H$  we multiply the equation  $\partial_t^{\alpha}u+\Delta^2 u=0 $ by $e_n$ and take into account that the operator $A$ is self-adjoint, so we have
\begin{equation*}
0=\langle\partial_t^{\alpha}u, e_n\rangle+\langle A u, e_n\rangle
=\partial_t^{\alpha}u_n+\langle  u, A e_n\rangle
=\partial_t^{\alpha}u_n+\lambda_n\langle  u, e_n\rangle
=\partial_t^{\alpha}u_n+\lambda_n u_n\,.
\end{equation*}
Therefore, in virtue of the initial conditions
\eqref{eq:ini-cond} $u_n(t)$ is the solution of  problem
\begin{equation}\label{Aeq:cauchy10}
\begin{cases}
\displaystyle
\partial_t^{\alpha}u_n(t)+\lambda_n u_n(t)=0\,,
\quad t\ge0,
\\
u_n(0)=\langle u_0,e_n\rangle,\quad
u_n'(0)=\langle u_1,e_n\rangle,
\end{cases}
\end{equation}
and hence by Lemma \ref{le:x-y} we get
\begin{equation}\label{eq:u_n}
u_n(t)=\langle u_0,e_n\rangle\ E_{\alpha,1}(-\lambda_n t^\alpha) 
+\langle u_1,e_n\rangle\ t E_{\alpha,2}(-\lambda_n t^\alpha) ,
\qquad t\ge0
\,.
\end{equation}
(i) Now, we take $u_0\in D(\sqrt{A})$, $u_1\in H$ and show that 
\begin{equation*}
u(t)=\sum_{n=1}^\infty u_n(t)e_n,
\end{equation*}
with $u_n(t)$ given by \eqref{eq:u_n}, is a weak solution of \eqref{eq:weakp} satisfying the initial conditions \eqref{eq:ini-cond}. First, we note that for any $t\in [0,T]$ we have $ u(t)\in D(\sqrt{A})$. Indeed, since
\begin{equation*}
\|\sqrt{A}  u(t)\|^2
=\sum_{n=1}^\infty \lambda_n | u_n(t)|^2
\le
2\sum_{n=1}^\infty \lambda_n \big| \langle u_0,e_n\rangle E_{\alpha,1}(-\lambda_nt^\alpha) \big|^2
+2\sum_{n=1}^\infty \lambda_n \big| \langle u_1,e_n\rangle t E_{\alpha,2}(-\lambda_nt^\alpha)\big|^2
,
\end{equation*}
thanks to \eqref{eq:stimeE} we have
\begin{equation*}
\begin{split}
 \lambda_n \big| \langle u_0,e_n\rangle E_{\alpha,1}(-\lambda_nt^\alpha) \big|^2
 &\le C\lambda_n \big| \langle u_0,e_n\rangle  \big|^2,
\\
 \lambda_n \big|\langle u_1,e_n\rangle t E_{\alpha,2}(-\lambda_nt^\alpha)\big|^2
& \le C t^{2-\alpha} \big|\langle u_1,e_n\rangle \big|^2\frac{\lambda_nt^\alpha}{(1+\lambda_nt^\alpha)^2}
\le C t^{2-\alpha} \big|\langle u_1,e_n\rangle \big|^2,
 \end{split}
\end{equation*}
and hence, being $\alpha<2$, we get
\begin{equation}\label{eq:deltau}
\|\sqrt{A} u(t)\|^2
\le C\|\sqrt{A} u_0\|^2+C T^{2-\alpha} \| u_1\|^2.
\end{equation}
Following the same reasoning pursued to get \eqref{eq:deltau}, we obtain for any $n\in\N$ 
\begin{equation*}
\Big\|\sqrt{A} \sum_{k=n}^\infty u_k(t)e_k\Big\|^2
\le C\sum_{k=n}^\infty \lambda_k \big| \langle u_0,e_k\rangle  \big|^2+C T^{2-\alpha} \sum_{k=n}^\infty \big|\langle u_1,e_k\rangle \big|^2,
\end{equation*}
and hence
\begin{equation*}
\lim_{n\to\infty}\sup_{t\in [0,T]}\Big\|\sqrt{A} \sum_{k=n}^\infty u_k(t)e_k\Big\|=0\,.
\end{equation*}
As a consequence, the series 
$\sum_{n=1}^\infty u_n(t)e_n$ 
is convergent in $D(\sqrt{A})$ uniformly in $t\in [0,T]$, so $u\in C([0,T];D(\sqrt{A}))$.
Moreover, $u(0)=\sum_{n=1}^\infty \langle u_0,e_n\rangle e_n=u_0$.

Concerning formula \eqref{eq:def-u_t} for $u'$, thanks to \eqref{eq:Ea1} and \eqref{eq:Eaa1} with $k=1$ we note that for $u_n(t)$ given by \eqref{eq:u_n} we have 
\begin{equation}\label{eq:un'}
u'_n(t)
=
-\lambda_n \langle u_0,e_n\rangle t^{\alpha-1}E_{\alpha,\alpha}(-\lambda_nt^\alpha)
+\langle u_1,e_n\rangle E_{\alpha,1}(-\lambda_nt^\alpha).
\end{equation}
We prove that $u'$ is given by \eqref{eq:def-u_t} and belongs to
$ C([0,T];D(A^{-\theta}))$
for $\theta\in\big(\frac{2-\alpha}{2\alpha},\frac12\big)$. Indeed, if $u'_n(t)$ is given by \eqref{eq:un'}, then we have
\begin{equation*}
\Big\|\sum_{n=1}^\infty u_n'(t)e_n\Big\|_{D(A^{-\theta})}^2
=\sum_{n=1}^\infty\lambda_n^{-2\theta}\big|u_n'(t)\big|^2
\,.
\end{equation*}
Since $0<\theta<\frac12$, thanks to \eqref{eq:maxbeta} we have
\begin{multline*}
\lambda_n^{-2\theta}\big|-\lambda_n \langle u_0,e_n\rangle t^{\alpha-1}E_{\alpha,\alpha}(-\lambda_nt^\alpha)\big|^2
\\
\le Ct^{2\alpha\theta+\alpha-2} \lambda_n |\langle u_0,e_n\rangle \big|^2\bigg(\frac{(\lambda_n t^{\alpha})^{\frac12-\theta}}{1+\lambda_nt^\alpha}\bigg)^2
\le Ct^{2\alpha\theta+\alpha-2} \lambda_n |\langle u_0,e_n\rangle \big|^2.
\end{multline*}
Therefore, taking into account that $\theta>\frac{2-\alpha}{2\alpha}$, for any $t\in [0,T]$ we get
\begin{equation*}
\Big\|\sum_{n=1}^\infty u_n'(t)e_n\Big\|_{D(A^{-\theta})}^2
\le
CT^{2\alpha\theta+\alpha-2}\|\sqrt{A} u_0\|^2+C\|u_1\|^2
\,.
\end{equation*}
Arguing as above, we can show that
the series $\sum_{n=1}^\infty u_n'(t)e_n$
is convergent in $D(A^{-\theta})$ uniformly in $t\in [0,T]$.
For that reason, the function $u$ is differentiable and 
\begin{equation*}
u'(t)=\sum_{n=1}^\infty u_n'(t)e_n,
\end{equation*}
that is, formula \eqref{eq:def-u_t} holds for $u'$. Moreover, $u'\in C([0,T];D(A^{-\theta}))$ and 
\begin{equation*}
u'(0)=\sum_{n=1}^\infty \langle u_1,e_n\rangle e_n=u_1. 
\end{equation*}
In addition, 
\begin{multline}\label{eq:utL2}
\|u'(t)\|^2=\sum_{n=1}^\infty\big|u_n'(t)\big|^2
\\
\le
2\sum_{n=1}^\infty\big|-\lambda_n \langle u_0,e_n\rangle t^{\alpha-1}E_{\alpha,\alpha}(-\lambda_nt^\alpha)\big|^2
+2\sum_{n=1}^\infty\big|\langle u_1,e_n\rangle E_{\alpha,1}(-\lambda_nt^\alpha)\big|^2
\,.
\end{multline}
Using \eqref{eq:stimeE} and \eqref{eq:maxbeta} we obtain
\begin{equation*}
\big|-\lambda_n \langle u_0,e_n\rangle t^{\alpha-1}E_{\alpha,\alpha}(-\lambda_nt^\alpha)\big|^2
\le Ct^{\alpha-2} \lambda_n |\langle u_0,e_n\rangle \big|^2\bigg(\frac{(\lambda_n t^{\alpha})^{\frac12}}{1+\lambda_nt^\alpha}\bigg)^2
\le Ct^{\alpha-2} \lambda_n |\langle u_0,e_n\rangle \big|^2.
\end{equation*}
Putting the above estimate into \eqref{eq:utL2} we get
\begin{equation*}
\int_0^T\|u'(t)\|^2\ dt
\le
C\frac{T^{\alpha-1}}{\alpha-1}\|\sqrt{A} u_0\|^2+C\|u_1\|^2
\,,
\end{equation*}
so $u'\in L^2(0,T;H)$. 

To evaluate $I^{2-\alpha}\big(u'-u_1\big)(t)$,
where
\begin{equation*}
u'(t)-u_1
=\sum_{n=1}^\infty\big[-\lambda_n \langle u_0,e_n\rangle t^{\alpha-1}E_{\alpha,\alpha}(-\lambda_nt^\alpha)
+\langle u_1,e_n\rangle \big(E_{\alpha,1}(-\lambda_nt^\alpha)-1\big)\big]e_n
\,,
\end{equation*} 
we observe that
\begin{equation*}
{\cal L}\big[I^{2-\alpha}\big(t^{\alpha-1}E_{\alpha,\alpha}(-\lambda_nt^\alpha)\big)\big](z)=\frac{z^{\alpha-2}}{z^{\alpha}+\lambda_n}
={\cal L}\big[tE_{\alpha,2}(-\lambda_nt^\alpha)\big)\big](z),
\end{equation*}
\begin{multline*}
{\cal L}\big[I^{2-\alpha}\big(E_{\alpha,1}(-\lambda_nt^\alpha)-1\big)\big](z)=z^{\alpha-2}\Big(\frac{z^{\alpha-1}}{z^{\alpha}+\lambda_n}-\frac1z\Big)
\\
=-\lambda_n\frac{z^{\alpha-3}}{z^{\alpha}+\lambda_n}
=-\lambda_n{\cal L}\big[t^{2}E_{\alpha,3}(-\lambda_nt^\alpha)\big)\big](z).
\end{multline*}
Therefore, by Lerch's Theorem that assures the uniqueness of inverse Laplace transforms and repeating  the previous argumentations one  proves 
\begin{equation}\label{eq:I2alpha}
I^{2-\alpha}\big(u'-u_1\big)(t)=
-\sum_{n=1}^\infty\lambda_n\big[ \langle u_0,e_n\rangle t E_{\alpha,2}(-\lambda_nt^\alpha)
+ \langle u_1,e_n\rangle t^2 E_{\alpha,3}(-\lambda_nt^\alpha)\big]e_n
\,,
\end{equation}
and $I^{2-\alpha}\big(u'-u_1\big)(t)\in C([0,T];H)$.

Next, if $v=\sum_{n=1}^\infty\ \langle v,e_n\rangle e_n$ belongs to $D(\sqrt A)$, then we have
\begin{equation*}
\langle I^{2-\alpha}\big(u'-u_1\big)(t),v\rangle
=
-\sum_{n=1}^\infty\lambda_n^{\frac12} \big[\langle u_0,e_n\rangle t E_{\alpha,2}(-\lambda_nt^\alpha)
+ \langle u_1,e_n\rangle t^2 E_{\alpha,3}(-\lambda_nt^\alpha)\big] \lambda_n^{\frac12}\langle v,e_n\rangle
\end{equation*}
We observe that by \eqref{eq:Eaa1} for $k=1$ and $k=2$ we have
\begin{multline*}
\sum_{n=1}^\infty\lambda_n^{\frac12} \frac{d}{dt}\Big[\langle u_0,e_n\rangle t E_{\alpha,2}(-\lambda_nt^\alpha)
+\langle u_1,e_n\rangle t^2 E_{\alpha,3}(-\lambda_nt^\alpha)\Big] \lambda_n^{\frac12}\langle v,e_n\rangle
\\
=\sum_{n=1}^\infty\lambda_n^{\frac12}\big[ \langle u_0,e_n\rangle  E_{\alpha,1}(-\lambda_nt^\alpha)
+ \langle u_1,e_n\rangle t E_{\alpha,2}(-\lambda_nt^\alpha)\big] \lambda_n^{\frac12}\langle v,e_n\rangle
\,.
\end{multline*}
Thanks to  \eqref{eq:stimeE} and \eqref{eq:maxbeta} we get
\begin{multline*}
\sum_{n=1}^\infty\lambda_n\big| \langle u_0,e_n\rangle  E_{\alpha,1}(-\lambda_nt^\alpha)
+ \langle u_1,e_n\rangle t E_{\alpha,2}(-\lambda_nt^\alpha)\big|^2 
\\
\le
C\sum_{n=1}^\infty\lambda_n \big|\langle u_0,e_n\rangle\big|^2\frac1{ (1+\lambda_nt^\alpha)^2 }
+Ct^{2-\alpha}\sum_{n=1}^\infty\big| \langle u_1,e_n\rangle\big|^2\Big(\frac{ \lambda_n^{\frac12}t^{\frac\alpha2}}{ 1+\lambda_nt^\alpha}\Big)^2 
\\
\le
C\big(\|\sqrt A u_0\|^2+t^{2-\alpha}\|u_1\|^2\big)
\,.
\end{multline*}
Therefore, since $\alpha<2$ for any $t\in [0,T]$ we have
\begin{multline*}
\sum_{n=1}^\infty\Big|\lambda_n^{\frac12}\frac{d}{dt}\Big[ \langle u_0,e_n\rangle t E_{\alpha,2}(-\lambda_nt^\alpha)
+ \langle u_1,e_n\rangle t^2 E_{\alpha,3}(-\lambda_nt^\alpha)\Big] \lambda_n^{\frac12}\langle v,e_n\rangle\Big|
\\
\le
C\big(\|\sqrt A u_0\|^2+\|u_1\|^2+\|\sqrt A v\|^2\big)\,,
\end{multline*}
whence
\begin{multline*}
\frac{d}{dt}\langle I^{2-\alpha}\big(u'-u_1\big)(t),v\rangle
\\
=-\sum_{n=1}^\infty\lambda_n^{\frac12}\frac{d}{dt}\Big[ \langle u_0,e_n\rangle t E_{\alpha,2}(-\lambda_nt^\alpha)
+ \langle u_1,e_n\rangle t^2 E_{\alpha,3}(-\lambda_nt^\alpha)\Big] \lambda_n^{\frac12}\langle v,e_n\rangle
\\
=-\sum_{n=1}^\infty\lambda_n^{\frac12}\big[ \langle u_0,e_n\rangle  E_{\alpha,1}(-\lambda_nt^\alpha)
+ \langle u_1,e_n\rangle t E_{\alpha,2}(-\lambda_nt^\alpha)\big] \lambda_n^{\frac12}\langle v,e_n\rangle
\,.
\end{multline*}
On the other hand
\begin{equation*}
\langle\sqrt A u(t),\sqrt A v\rangle
=\sum_{n=1}^\infty\lambda_n^{\frac12}\big[ \langle u_0,e_n\rangle  E_{\alpha,1}(-\lambda_nt^\alpha)
+ \langle u_1,e_n\rangle t E_{\alpha,2}(-\lambda_nt^\alpha)\big] \lambda_n^{\frac12}\langle v,e_n\rangle
\,,
\end{equation*}
and hence for any $t\in [0,T]$ we have
\begin{equation*}
\frac{d}{dt}\langle I^{2-\alpha}\big(u'-u_1\big)(t),v\rangle+\langle\sqrt A u(t),\sqrt A v\rangle=0
\,,
\end{equation*}
that is \eqref{eq:w-int} holds. 

In conclusion, $u$ given by \eqref{eq:def-u0} is the weak solution of \eqref{eq:weakp} satisfying
 the initial conditions \eqref{eq:ini-cond}.

(ii) We assume that $u_0\in D({A})$ and $u_1\in D(\sqrt{A})$.
Repeating  similar argumentations  to those done before, involving again \eqref{eq:stimeE} and \eqref{eq:maxbeta}, we  have that the function $u$ given by
\eqref{eq:def-u0} belongs to $C([0,T];D(A))\cap C^1([0,T];H)$.

From \eqref{eq:I2alpha} it follows that
\begin{equation}\label{eq:dtI2}
\frac{d}{dt} I^{2-\alpha}\big(u'-u_1\big)(t)=-\sum_{n=1}^\infty\lambda_n\big[ \langle u_0,e_n\rangle E_{\alpha,1}(-\lambda_nt^\alpha)
+\langle u_1,e_n\rangle t E_{\alpha,2}(-\lambda_nt^\alpha)\big]e_n
\end{equation}
and $\frac{d}{dt} I^{2-\alpha}\big(u'-u_1\big)(t)\in C([0,T];H)$. Indeed, for any $t\in[0,T]$
\begin{equation*}
\lambda_n^2\big| \langle u_0,e_n\rangle E_{\alpha,1}(-\lambda_nt^\alpha)
+\langle u_1,e_n\rangle t E_{\alpha,2}(-\lambda_nt^\alpha)\big|^2
\le
C\big(\lambda_n^2 |\langle u_0,e_n\rangle|^2+T^{2-\alpha}\lambda_n |\langle u_1,e_n\rangle|^2\big)
\,,
\end{equation*}
whence
\begin{equation*}
\sum_{n=1}^\infty\big|\lambda_n \langle u_0,e_n\rangle E_{\alpha,1}(-\lambda_nt^\alpha)
+\lambda_n\langle u_1,e_n\rangle t E_{\alpha,2}(-\lambda_nt^\alpha)\big|^2
\le
C\big(\|A u_0\|^2+T^{2-\alpha}\|\sqrt Au_1\|^2\big)
\,.
\end{equation*}
From \eqref{eq:un'}, taking into account \eqref{eq:Eaaa} and \eqref{eq:Ea1}, we have
\begin{equation}\label{eq:un''}
u''_n(t)
=
-\lambda_n\big[ \langle u_0,e_n\rangle t^{\alpha-2}E_{\alpha,\alpha-1}(-\lambda_nt^\alpha)
+\langle u_1,e_n\rangle t^{\alpha-1} E_{\alpha,\alpha}(-\lambda_nt^\alpha)\big].
\end{equation}
Moreover,
\begin{equation}
\big\|\sum_{n=1}^\infty u''_n(t)e_n\big\|^2
\le
C(t^{2\alpha-4}\|A u_0\|^2+t^{\alpha-2}\|\sqrt A u_1\|^2),
\end{equation}
and hence
\begin{equation}
\big\|\sum_{n=1}^\infty u''_n(t)e_n\big\|
\le
C(t^{\alpha-2}\|A u_0\|+t^{\frac\alpha2-1}\|\sqrt A u_1\|).
\end{equation}
Therefore,
\begin{equation}
\int_0^T\big\|\sum_{n=1}^\infty u''_n(t)e_n\big\|\ dt
\le
C(T^{\alpha-1}\|A u_0\|+T^{\frac\alpha2}\|\sqrt A u_1\|),
\end{equation}
that is $\sum_{n=1}^\infty u''_n(t)e_n$ belongs to $L^1(0,T;H)$.
Therefore, for any $t\in [0,T]$ we have
\begin{equation}
\int_0^t\sum_{n=1}^\infty u''_n(s)e_n\ ds=\sum_{n=1}^\infty\int_0^t u''_n(s) ds \ e_n
=\sum_{n=1}^\infty( u'_n(t) -\langle u_1,e_n\rangle) \ e_n=u'(t)-u_1.
\end{equation}
Since $u'$ is absolutely continuous, taking into account \eqref{eq:der-fracI1}, we have
\begin{equation}\label{eq:abs-c}
\frac{d}{dt} I^{2-\alpha}\big(u'-u_1\big)(t)=\partial_t^{\alpha}u(t)
\,.
\end{equation}
In particular, $\partial_t^{\alpha}u\in C([0,T];H)$ and, thanks to \eqref{eq:dtI2},  formula \eqref{eq:def-u-alpha} holds.
%
%
%
Finally, from \eqref{eq:w-int} and \eqref{eq:abs-c} we have
for any $v\in D(\sqrt{A})$ and $t\in (0,T)$
\begin{equation*}
0=\langle \frac{d}{dt}I^{2-\alpha}\big(u'-u_1\big)(t),v\rangle+\langle \sqrt{A}u(t),\sqrt{A} v \rangle=
\langle  \partial_t^{\alpha}u(t),v\rangle+\langle Au(t), v \rangle=\langle  \partial_t^{\alpha}u(t)+Au(t),v\rangle,
\end{equation*}
that is $\partial_t^{\alpha}u+A u=0$ in 
$(0,T)$. 

In conclusion, $u$ satisfies the conditions of Definition \ref{de:wss}-2, that is, $u$ is a strong solution of \eqref{eq:weakp}.
\end{Proof}

\section{Examples}
In this section we give an application of our well-posedness results by discussing two examples of concrete models involving well-known  partial differential operators. Throughout the section, we denote by $\Omega$ a bounded open domain in $\R^N$, $N\ge 1$, with sufficiently smooth boundary $\partial\Omega$.
In both examples we consider the Hilbert space  $H=L^2(\Omega)$,
endowed with the inner product and norm defined by
\begin{equation*}
\langle u,v\rangle=\int_{\Omega}u(x)v(x)\ dx,
\qquad
\|u\|=\left(\int_{\Omega}|u(x)|^{2}\ dx\right)^{1/2}\qquad
u,v\in L^2(\Omega).
\end{equation*}
\subsection{Time-fractional wave equations}
We analyse the fractional boundary value problem
\begin{equation}\label{eq:weakwE}
\begin{cases}
\partial_t^{\alpha}u=\Delta u
\qquad \mbox{in}\ (0,T)\times\Omega,
\\
u=0    \hskip1.6cm \mbox{on}\ (0,T)\times\partial\Omega.
\end{cases}
\end{equation}
We  rewrite \eqref{eq:weakwE} as an abstract equation of the type \eqref{eq:weakpI} by introducing the operator $A$ as follows:
\begin{equation*}
\begin{split}
D(A)&= H^2(\Omega)\cap H^1_0(\Omega)
\\
(Au)(x)&=-\Delta u(x), 
\quad u\in D(A),
\  x\in\Omega.
\end{split}
\end{equation*}
It is well known that $A$  is a linear self-adjoint positive
  operator on $L^2(\Omega)$ with dense domain.
Moreover, the fractional power $\sqrt A$ of $A$ is well defined and $D(\sqrt A) = H^1_0(\Omega)$.

We assume that the eigenvalues $\lambda_n$, $n\in\N$, of the operator $A$ are all distinct numbers, whence the eigenspace generated by 
$\lambda_n$  has dimension one.

We are in conditions to apply Theorem \ref{th:reg-l2} to get a well-posedness result for problem \eqref{eq:weakwE}.

\begin{theorem}
\begin{enumerate}
\item[(i)] 
If $u_0\in H^1_0(\Omega)$ and $u_1\in L^{2}(\Omega)$, then the function
\begin{equation}\label{eq:def-u0E1}
u(t,x)=\sum_{n=1}^\infty\big[  \langle u_0,e_n\rangle E_{\alpha,1}(-\lambda_nt^\alpha) 
+\langle u_1,e_n\rangle t E_{\alpha,2}(-\lambda_nt^\alpha)\big]e_n(x)
\end{equation}
is the unique weak solution of \eqref{eq:weakwE} satisfying
 the initial conditions
\begin{equation*}
u(0,\cdot)=u_{0},\quad
u_t(0,\cdot)=u_{1}.
\end{equation*}
In addition
\begin{equation*}
u_t(t,x)=\sum_{n=1}^\infty\big[-\lambda_n \langle u_0,e_n\rangle t^{\alpha-1}E_{\alpha,\alpha}(-\lambda_nt^\alpha)
+\langle u_1,e_n\rangle E_{\alpha,1}(-\lambda_nt^\alpha)\big]e_n(x)
\,,
\end{equation*}
and $u_t\in C([0,T];D(A^{-\theta}))$ for  $\theta\in\big(\frac{2-\alpha}{2\alpha},\frac12\big)$.
\item[(ii)]
For $u_0\in H^2(\Omega)\cap H^1_0(\Omega)$ and $u_1\in H^1_0(\Omega)$
the weak solution given by \eqref{eq:def-u0E1} is a strong one and 
\begin{equation*}
\partial_t^{\alpha}u(t,x)=-\sum_{n=1}^\infty\big[\lambda_n \langle u_0,e_n\rangle E_{\alpha,1}(-\lambda_nt^\alpha)
+\lambda_n\langle u_1,e_n\rangle t E_{\alpha,2}(-\lambda_nt^\alpha)\big]e_n(x)
\,.
\end{equation*}
\end{enumerate}

\end{theorem}
\subsection{Time-fractional Petrovsky systems}
We consider the system
\begin{equation}\label{eq:weakpE}
\begin{cases}
\partial_t^{\alpha}u+\Delta^2 u=0
\qquad \mbox{in}\ (0,T)\times\Omega,
\\
u=\Delta u=0    \hskip1.4cm \mbox{on}\ (0,T)\times\partial\Omega.
\end{cases}
\end{equation}
We can recast \eqref{eq:weakpE} as an abstract problem by defining the operator $A$ in this way:
\begin{equation}\label{eq:defa}
\begin{split}
D(A)&= H^4(\Omega)\cap \{u\in H^3(\Omega): u=\Delta u=0\ \mbox{on}\ \partial\Omega\}
\\
(Au)(x)&=\Delta^2 u(x), 
\quad u\in D( A),
\  x\in\Omega.
\end{split}
\end{equation}
We assume that the eigenvalues $\lambda_n'$, $n\in\N$, of the operator $-\Delta$ with domain $H^2(\Omega)\cap H^1_0(\Omega)$ are all distinct numbers, and hence the eigenspace generated by 
$\lambda_n'$  has dimension one.

The biharmonic operator given by \eqref{eq:defa} satisfies: 

\begin{itemize}
\renewcommand\labelitemi{--}
\item $A$ is self-adjoint, because,
integrating by parts and taking into account the boundary conditions satisfied by the elements of $D(A)$, we have
\begin{equation*}
\langle A u,v\rangle=\int_{\Omega}\Delta u(x)\Delta v(x)\ dx=\langle u,A v\rangle
\qquad
u,v\in D(A);
\end{equation*}
\item
the fractional power $\sqrt A$ of the operator $A$ is $-\Delta$ with domain $D(\sqrt A)=H^2(\Omega)\cap H^1_0(\Omega)$;
\item $A$ is positive, since \eqref{eq:opA} is satisfied. Indeed,
\begin{equation*}
\|\Delta u\|\cdot\|u\|\ge\langle -\Delta u,u\rangle=\int_{\Omega}|\nabla u(x)|^2\ dx\ge C\int_{\Omega}|u(x)|^2\ dx
\end{equation*}
and hence
\begin{equation*}
\langle Au,u\rangle=\int_{\Omega}|\Delta u(x)|^2\ dx=\|\Delta u\|^2\ge C\| u\|^2
\qquad
u\in D(A);
\end{equation*}
\item the domain $D(A)$ is dense in $L^2(\Omega)$  by the density of $C^\infty_c(\Omega)$ in $L^2(\Omega)$.

\end{itemize}

\begin{lemma}
The spectrum of the operator $A$ consists of the sequence tending to $+\infty$ of eigenvalues 
$\lambda_n=(\lambda'_n)^2$, where $\lambda'_n$  are the eigenvalues of the operator
$-\Delta$ in $H^2(\Omega)\cap H^1_0(\Omega)$.

The eigenfunctions $e_n$ of $-\Delta$  \big($-\Delta e_n=\lambda'_n e_n$\big), which constitutes an orthonormal basis of $L^2(\Omega)$,
are also eigenfunctions of $\Delta^2$ \big($\Delta^2e_n=\lambda_n e_n$\big).
\end{lemma}
\begin{Proof}
It is enough to note that
\begin{equation*}
\Delta^2e_n= -\Delta(-\Delta e_n)= -\lambda_n'\Delta(e_n)=(\lambda_n')^2 e_n
\,.
\end{equation*}
The operator $\Delta^2$ cannot have other eigenfunctions, because the sequence $\{e_n\}$  of the eigenfunctions of $-\Delta$  constitutes an orthonormal basis of $L^2(\Omega)$.
\end{Proof}

In the case of the biharmonic operator a weak solution can be called a $H^2$-solution. This terminology is suggested by the analysis of the stationary case given in  \cite {NS}. 

Finally, by Theorem \ref{th:reg-l2} a well-posedness result for system \eqref{eq:weakpE} follows.

\begin{theorem}
\begin{enumerate}
\item[(i)] 
If $u_0\in H^2(\Omega)\cap H^1_0(\Omega)$ and $u_1\in L^{2}(\Omega)$, then the function
\begin{equation}\label{eq:def-u0E}
u(t,x)=\sum_{n=1}^\infty\big[  \langle u_0,e_n\rangle E_{\alpha,1}(-\lambda_nt^\alpha) 
+\langle u_1,e_n\rangle t E_{\alpha,2}(-\lambda_nt^\alpha)\big]e_n(x)
\end{equation}
is the unique $H^2$-solution of \eqref{eq:weakpE} satisfying
 the initial conditions
\begin{equation*}
u(0,\cdot)=u_{0},\quad
u_t(0,\cdot)=u_{1}.
\end{equation*}
In addition
\begin{equation*}
u_t(t,x)=\sum_{n=1}^\infty\big[-\lambda_n \langle u_0,e_n\rangle t^{\alpha-1}E_{\alpha,\alpha}(-\lambda_nt^\alpha)
+\langle u_1,e_n\rangle E_{\alpha,1}(-\lambda_nt^\alpha)\big]e_n(x)
\,,
\end{equation*}
and $u_t\in C([0,T];D(A^{-\theta}))$ for  $\theta\in\big(\frac{2-\alpha}{2\alpha},\frac12\big)$.
\item[(ii)]
For $u_0\in H^4(\Omega)\cap \{u\in H^3(\Omega): u=\Delta u=0\ \mbox{on}\ \partial\Omega\}$ and $u_1\in H^2(\Omega)\cap H^1_0(\Omega)$
the $H^2$-solution given by \eqref{eq:def-u0E} is a strong one and 
\begin{equation*}
\partial_t^{\alpha}u(t,x)=-\sum_{n=1}^\infty\big[\lambda_n \langle u_0,e_n\rangle E_{\alpha,1}(-\lambda_nt^\alpha)
+\lambda_n\langle u_1,e_n\rangle t E_{\alpha,2}(-\lambda_nt^\alpha)\big]e_n(x)
\,.
\end{equation*}
\end{enumerate}

\end{theorem}


\begin{thebibliography}{99}

\itemsep=\smallskipamount






\bibitem{B} H. Brezis, {\em Analyse fonctionnelle. Th\'eorie et applications.}  Collection Math\'ematiques Appliqu\'ees pour la Ma\^\i trise. Masson, Paris, 1983. 








\bibitem{Fujita} Y. Fujita, Integrodifferential equation which interpolates the heat equation and the wave equation, I, II, Osaka J. Math. 27 (1990) 309--321, 797--804.




\bibitem{GLY} 
R. Gorenflo, Y. Luchko, M. Yamamoto, {\em Time-fractional diffusion equation in the fractional Sobolev spaces} Fract. Calc. Appl. Anal. 18 (2015), no. 3, 799--820.







\bibitem{LSw} P. Loreti, D. Sforza,  {\em Fractional diffusion-wave equation: hidden regularity for weak solutions}
arXiv:2007.13581 

\bibitem{HMP} H. R. Henr\'iquez, J. G. Mesquita, J. C. Pozo, 
Existence of solutions of the abstract Cauchy problem of fractional order.
J. Funct. Anal. 281 (2021), no. 4, 109028. 

%


\bibitem{Lunardi} A. Lunardi,  Interpolation theory. Third edition [of MR2523200]. Appunti. Scuola Normale Superiore di Pisa (Nuova Serie) [Lecture Notes. Scuola Normale Superiore di Pisa (New Series)], 16. Edizioni della Normale, Pisa, 2018. xiv+199 pp. 


\bibitem{MG} F. Mainardi, R. Gorenflo,
{\em Time-fractional derivatives in relaxation processes: a tutorial survey} Fract. Calc. Appl. Anal. 10 (2007), no. 3, 269--308.

\bibitem{NS} S. A. Nazarov, G. Sweers, A hinged plate equation and iterated Dirichlet Laplace operator on domains with concave corners, J. Differential Equations 233 (2007), no. 1, 151--180. 



\bibitem{Pazy} A. Pazy, Semigroups of Linear Operators and Applications to Partial Differential Equations, Springer-Verlag, Berlin, 1983.

\bibitem{Pod} I. Podlubny, Fractional Differential Equations, Academic Press, San Diego, 1999.



\bibitem{SY}
K. Sakamoto, M. Yamamoto, {\em Initial value/boundary value problems for fractional diffusion-wave equations and applications to some inverse
problems}, J. Math. Anal. Appl. 382 (1) (2011) 426--447.

\bibitem{KY}
Y. Kian, M. Yamamoto, Well-posedness for weak and strong solutions of non-homogeneous initial boundary value problems for fractional diffusion equations. Fract. Calc. Appl. Anal. 24 (2021), no. 1, 168--201.


\end{thebibliography}
\end{document}